\documentclass[11pt]{article}

\usepackage[top=4cm, bottom=2cm, left=2.2cm, right=2.2cm]{geometry}
\usepackage{subfig}  
\usepackage{amssymb}
\usepackage{graphicx}
\usepackage{amsmath}
\usepackage[titletoc]{appendix}
\usepackage{appendix}
\usepackage{pifont}
\usepackage{bbm}
\usepackage{amsfonts}
\usepackage{mathrsfs}
\usepackage{amssymb}
\usepackage{amsmath}
\usepackage{a4wide}
\usepackage{stmaryrd}
\usepackage{txfonts}
\usepackage{tabularx}
\usepackage{booktabs}
\usepackage{ctable}

\usepackage{latexsym}

\let\er\eqref

\newcommand{\R}{{\mathbb R}}

\newtheorem{thm}{Theorem}
\newtheorem{lemma}[thm]{Lemma}
\newtheorem{proposition}[thm]{Proposition}
\newtheorem{remark}[thm]{Remark}

\newtheorem{definition}[thm]{Definition}

\begin{document}
\fontfamily{ptm}\fontseries{sb}\selectfont
\title{
{\bf\Large Chemotaxis model with subcritical exponent in nonlocal reaction}
 \\
}
\author{ Shen Bian \footnote{Beijing University of Chemical Technology, 100029, Beijing. Universit\"at Mannheim, 68131, Mannheim. Email: \texttt{bianshen66@163.com}. Partially supported by National Science Foundation of China (Grant No. 11501025) and the Alexander von Humboldt Foundation.}
  \and Li Chen\footnote{Universit\"at Mannheim, 68131, Mannheim. Email: \texttt{chen@math.uni-mannheim.de}. Partially supported by DFG Project CH 955/3-1.}
  \and Evangelos A. Latos\footnote{Universit\"at Mannheim, 68131, Mannheim. Email: \texttt{evangelos.latos@math.uni-mannheim.de}. Partially supported by DFG Project CH 955/3-1.}
}
\date{}
\maketitle
\begin{center}
{\bf\small Abstract}

\vspace{3mm}\hspace{.05in}\parbox{5in} {This paper deals with a parabolic-elliptic chemotaxis system with nonlocal type of source in the whole space. It's proved that the initial value problem possesses a unique global solution which is uniformly bounded. Here we identify the exponents regimes of nonlinear reaction and aggregation in such a way that their scaling and the diffusion term coincide (see Introduction). Comparing to the classical KS model (without the source term), it's shown that how energy estimates give natural conditions on the nonlinearities implying the absence of blow-up for the solution without any restriction on the initial data.
}
\end{center}


\section{Introduction}
\def\theequation{1.\arabic{equation}}\makeatother
\setcounter{equation}{0}
\def\thetheorem{1.\arabic{thm}}\makeatother
\setcounter{thm}{0}

In this work, we analyze qualitative properties of non-negative solutions for the chemotaxis system in dimension $n \ge 3$ with linear diffusion given by
\begin{align}\label{star00}
\left\{
  \begin{array}{ll}
 u_t=\Delta u-\nabla \cdot (u^{\sigma} \nabla v)+u^\alpha \left( 1-\int_{\R^n} u^\beta dx \right),\quad   &x \in \R^n, t>0,
\\
 u(x,0)=u_0(x)\geq 0,\quad & x \in \R^n.
  \end{array}
\right.
\end{align}
Here $v(x,t)$ expresses the chemical substance concentration and it is given by the fundamental solution
\begin{align}
v(x,t)=K \ast u(x,t) =c_n \int_{\R^n} \frac{u(y,t)}{|x-y|^{n-2}}dy
\end{align}
where
\begin{align}
c_n=\frac{1}{n(n-2)b_n},~~b_n=\frac{\pi^{n/2}}{\Gamma(n/2+1)},
\end{align}
$b_n$ is the volume of $n$-dimensional unit ball. This system without the reaction has been proposed as a model for chemotaxis driven cell movement \cite{KS70}. Here $\sigma \ge 1$ in chemotaxis is to model the nonlinear aggregation, $u^\alpha \left( 1-\int_{\R^n} u^\beta dx \right)$ with $\alpha>1, \beta>1$ is the reaction term representing nonlinear growth under nonlocal resource consumption of the bacteria \cite{F37,KPP37}.

Initial data will be assumed
\begin{align}\label{initial}
u_0 \ge 0,~~ u_0 \in W^{2,n-1}(\R^n) \cap L^1(\R^n).
\end{align}
Under the initial assumptions, we consider the case $\sigma \ge 1$ and prove that in the following cases,
\begin{itemize}
  \item[(1)] $\sigma+1 \le \alpha$ and $\alpha<1+2\beta/n$,
  \item[(2)] $\sigma+1>\alpha$ and $(\sigma+1)(n+2)<2 \beta+2 \alpha+n$,
\end{itemize}
the solution of \er{star00} is unique and global without any restriction on the initial data.

In the follows, we define the exponent $p$ arising from Sobolev inequality \cite{Lieb202} and the notation $Q_T$
\begin{align}
p&:=\frac{2n}{n-2},~~\label{p}\\
Q_T&:=\R^n \times (0,T) \mbox{~~for~all~} T>0. \nonumber
\end{align}
Throughout this paper, we deal with a strong solution of \er{star00} which is defined:
\begin{definition}\label{def1}
Let $\sigma \ge 1, \alpha>1, \beta>1$, $u_0$ satisfies \er{initial}. $u(x,t)$ is called a strong solution on [0,T) if
\begin{itemize}
  \item[(1)] $u \in W^{2,n-1}(Q_T),~u_t\in L^{n-1}(Q_T)$,
  \item[(2)] $u \in L^\infty(0,T;L^1\cap L^\infty(\R^n))$.
\end{itemize}
\end{definition}
Our main result concerning the solution can be summarized as follows:
\begin{thm}\label{thm1}
Let $\alpha>1, \beta>1, \sigma \ge 1$. If either
$$\sigma+1 \le \alpha <1+\frac{2 \beta}{n}$$ or $$\alpha<\sigma+1<\frac{2(\beta+\alpha)}{n+2}+\frac{n}{n+2},$$ then for any initial data satisfying \er{initial}, problem \er{star00} possesses a unique global strong solution defined by Definition \ref{def1} which is uniformly bounded, i.e., for any $t>0,$ then
\begin{align}
\|u(\cdot,t)\|_{L^q(\R^n)} \le C_1(\|u_0\|_{L^q(\R^n)}),~~q \in [\beta+\alpha-1,\infty], \\
\|u(\cdot,t)\|_{L^q(\R^n)} \le C_2(\|u_0\|_{L^q(\R^n)},t),~~q \in [1,\beta+\alpha-1).
\end{align}
Here $C_1$ is a positive constant only depending on $\|u_0\|_{L^q(\R^n)}$ but not on $t$.
\end{thm}

Actually, problem \er{star00} contains three terms, the diffusion term $\Delta u$, the nonlocal aggregation term $\nabla\cdot(u^\sigma \nabla v)$ and the nonlinear growth term $u^\alpha (1-\int_{\R^n} u^\beta dx)$ (where $-u^\alpha \int_{\R^n} u^\beta dx$ can be viewed as the death contributing to the global existence), then the competition arises between the diffusion, the death and the aggregation, the growth. Indeed, \er{star00} can be recast as
\begin{align}\label{starstar00}
u_t=\Delta u-\sigma u^{\sigma-1}\nabla u \cdot \nabla v+u^{\sigma+1}+u^\alpha\left(1-\int_{\R^n} u^\beta dx\right),
\end{align}
if $\sigma+1=\alpha,$ then the particular nonlinear reaction exponent
$$\alpha=2\beta/n+1$$
gives the balance of diffusion and aggregation, reaction. In fact, plugging $u_\lambda(x,t)=\lambda^{\frac{n}{\beta}}u(\lambda x, \lambda^2 t)$ into \er{star00}, it's easy to verify that $u_\lambda(x,t)$ is also a solution of \er{star00} and the scaling preserves the $L^\beta$ norm in space, the diffusion term $\lambda^{n/\beta+2}\Delta u(\lambda x, \lambda^2 t)$ has the same scaling as the aggregation term $\lambda^{(\sigma+1)n/\beta}\nabla \cdot (u \nabla (K\ast u))(\lambda x, \lambda^2 t)$ and the reaction term $\lambda^{n\alpha/\beta} u^\alpha (1-\int_{\R^n}u^\beta dx)(\lambda x, \lambda^2 t)$ if and only if $\alpha=2\beta/n+1$. From observing the rescaled equation we can see that when $n(\alpha-1)/\beta<2,$ for low density (small $\lambda$), the aggregation dominates the diffusion thus prevents spreading. While for high density (large $\lambda$), the diffusion dominates the aggregation and thus blow-up is precluded. Hence, in this case, the solution will exist globally (Theorem \ref{thm1}). On the other hand, if $n(\alpha-1)/\beta>2,$ then the diffusion dominates for low density and the density had infinite-time spreading, the aggregation manipulates for high density and the density has finite time blow-up. Therefore, our conjecture is that there exists finite time blow-up for $\alpha-1>2\beta/n.$ Moreover, for $\alpha-1=2\beta/n$, similar to \cite{Jose2009}, we guess that there is a critical value for the initial data sharply separating global existence and finite time blow-up.

Moreover, noticing that \er{starstar00} includes $u_t, \Delta u, u^{\sigma+1}, u^\alpha, \nabla \cdot (u^\sigma \nabla v)$. Therefore, we detect the following equations are similar to \er{starstar00}:
\begin{align}\label{DR}
\left\{
  \begin{array}{ll}
    u_t=\Delta u+u^\alpha, & x \in \R^n, t>0, \\
    u(x,0)=u_0(x), & x \in \R^n.
  \end{array}
\right.
\end{align}
and
\begin{align}\label{DA}
\left\{
  \begin{array}{ll}
    u_t=\Delta u-\nabla \cdot (u^\sigma \nabla v), & x \in \R^n, t>0, \\
    u(x,0)=u_0(x), & x \in \R^n.
  \end{array}
\right.
\end{align}
In order to compare \er{DR}, \er{DA} with \er{starstar00}, we take $\sigma+1=\alpha$ (for simplicity) to find the effects of the nonlocal reaction $u^\alpha \left(1-\int_{\R^n} u^\beta dx \right)$. Concerning \er{starstar00} and \er{DR}, in this paper, we prove that for $\alpha<1+2\beta/n,$ \er{star00} admits a unique and global solution. While for the Fujita type equation \er{DR}, it's known that for $\alpha<1+2/n,$ there is no global solution \cite{Fu66} (For comparison, $\lambda^{\frac{n}{\beta}}u(\lambda x, \lambda^2 t)$ in \er{star00} is just the mass invariant scaling). As to \er{DA} and \er{starstar00}, the most remarkable difference is that the mass conservation holds for \er{DA} but not for \er{starstar00}, using this property it's been proved that the solution of \er{DA} exists globally with small initial data \cite{BL13,BDP06,JL92,P07,SK06}, while \er{starstar00} has a unique and global solution without any restriction on the initial data. Thus we conclude that the reaction term can prevent blow-up.

In addition, Keller-Segel model with local reaction term in bounded domain has been widely studied by virtue of comparison principle and those estimates that are valid in bounded domain \cite{GST16,HZ16,IS16,NT13,TW07,WL17,WMZ14,ZhLi15,Zh15}, that's
\begin{align} \label{star11}
\left\{
      \begin{array}{ll}
       u_t=\Delta u-\chi \nabla \cdot (u^{\sigma} \nabla c)+f(u), &  x \in \Omega, t>0, \\
       -\Delta c+ c=u, & x \in \Omega, t>0, \\
       \nabla u \cdot \vec{n}=\nabla v \cdot \vec{n}=0, &  x \in \partial \Omega.
      \end{array}\right.
\end{align}
For $\sigma=1,$ \cite{TW07} showed that model \er{star11} with $f(u) \le a-b u^2, u \ge 0, a,b>0$ possesses a global bounded solution under the assumption $b>\frac{n-2}{n}.$ The case $\sigma>1$ was considered in \cite{GST16} with $f(u)=\mu u(1-u^\alpha)$, if $\alpha>\sigma$ or $\alpha=\sigma$ and $\mu>\frac{n\alpha -2}{n\alpha+2(\sigma-1)}$, then \er{star11} admits a unique global solution. For other similar cases, one can refer to \cite{HT17,NT13,WMZ14}.

In brief, comparing to the above models, the absence of mass conservation (model \er{DA}) and the comparison principle (model \er{DR},\er{star11}) are two obstacles in our model \er{star00} \cite{Fu66,P07,SK06,WMZ14}. Besides, the nonlocal reaction makes the key energy estimates more difficult and many tools in the bounded domain  can't work in the whole space \cite{GST16,NT13,TW07,WMZ14}. In our results, we will use analytical methods in the energy estimates and derive the conditions on $\alpha,\beta,\sigma$ for global existence (Theorem \ref{thm1}).

The main work is devoted to the global unique solution of model \er{star00} for $\alpha>1, \beta>1, \sigma \ge 1$, with that aim Section \ref{sec3} considers the local existence and uniqueness of the solution. In Section \ref{sec4}, the a priori estimates are performed and show that $-u^\alpha \int_{\R^n} u^\beta dx$ plays a crucial role on the global existence, thus complete the proof of Theorem \ref{thm1}. Here we split the arguments into several parts strongly depending on the exponents $\alpha,\beta$ and $\sigma,$ consequently the uniformly boundedness is obtained by virtue of the Moser iterative method. Section \ref{sec5} discusses some open questions of Eq. \er{star00}.

\section{Preliminaries} \label{sec2}
\def\theequation{2.\arabic{equation}}\makeatother
\setcounter{equation}{0}
\def\thetheorem{2.\arabic{thm}}\makeatother

We firstly state some lemmas which will be used in the proof of local existence and Theorem 2.

Consider the Cauchy problem
\begin{align}\label{Cauchy}
\left\{
  \begin{array}{ll}
    z_t=\Delta z-\nabla \cdot (z H)+F,~~x \in \R^n, t>0, \\
    z(x,0)=z_0(x).
  \end{array}
\right.
\end{align}
Then the solution $z(x,t)$ can be expressed from semigroup theory \cite{Pa83} as follows:
\begin{lemma}\label{lem21}
Let $X$ be a Banach space, $z_0 \in X, H \in L^\infty(0,T;X)$ and $F\in L^\infty(0,T;X)$. $z(x,t)$ is given by \cite{BDP06,SK06}
\begin{align}
z(x,t)=G(\cdot,t)\ast z_0+\int_0^t \nabla G(\cdot,t-s) \ast [z(\cdot,s)H(\cdot,s)] ds+\int_0^t G(\cdot,t-s)\ast F(\cdot,s)ds,~~0 \le t \le T
\end{align}
is the mild solution of \er{Cauchy} on $[0,T].$ Here $G(x,t)=\frac{1}{(4 \pi t)^{\frac{n}{2}}}e^{-\frac{|x|^2}{4t}}$ is the Green function associated to the heat equation.
\end{lemma}

The following lemma is an immediate consequence from Sobolev inequality \cite{Lieb202} which will play an important role in the proof of global existence of solutions for equation \er{star00}.
\begin{lemma}[\cite{BL16}]\label{interpolation}
Let $p$ is expressed by \er{p}, $1 \le r<q<p$ and $\frac{q}{r}<\frac{2}{r}+1-\frac{2}{p}$, then for any $v\in H^1(\R^n)$ and $v \in L^r(\R^n)$, it holds
\begin{align}
\|v\|_{L^q(\R^n)}^q \le C(n) C_0^{-\frac{ \lambda q}{2-\lambda q}} \|v\|_{L^r(\R^n)}^{\gamma} + C_0 \|\nabla v\|_{L^2(\R^n)}^2,~~n \ge 3, \label{nge3}
\end{align}
Here $C(n)$ is a constant depending on $n$, $C_0$ is an arbitrarily positive constant and
\begin{align}
\lambda=\frac{\frac{1}{r}-\frac{1}{q}}{\frac{1}{r}-\frac{1}{p}} \in (0,1),~~
\gamma=\frac{2(1-\lambda) q}{2-\lambda q}=\frac{2\left(1-\frac{q}{p}\right)}{\frac{2-q}{r}-\frac{2}{p}+1}.
\end{align}
\end{lemma}

\section{Local existence and uniqueness} \label{sec3}
\def\theequation{3.\arabic{equation}}\makeatother
\setcounter{equation}{0}
\def\thetheorem{3.\arabic{thm}}\makeatother

This part concerns local existence of the strong solution of \er{star00}. The result is standard, more detailed arguments can be found in \cite{BDP06,SK06,TW07}.
\begin{proposition}\label{pro23}
Let $\alpha>1, \beta>n/2, \sigma \ge 1$. Assume that the initial data $u_0\in W^{2,n-1}(\R^n) \cap L^1(\R^n),$ then there exists a maximal existence time $T_{\max} \in (0,\infty]$ such that $u(x,t) \in W^{2,n-1}(Q_T) \cap L^\infty(0,T;L^1(\R^n))$ is the unique non-negative strong solution of problem \er{star00}. Furthermore, if $T_{\max}<+\infty,$ then
\begin{align}\label{blowup}
\displaystyle \lim_{t \to T_{\max}} \|u(\cdot,t)\|_{L^\infty(\R^n)}=\infty.
\end{align}
\end{proposition}
\begin{remark}
Here $\beta$ is chosen to be $\beta>n/2$ in order to prove the local existence and the a priori estimates Proposition \ref{pro25}. By sobolev embedding theorem, $W^{2,n-1}(\R^n)$ embedding into $L^\infty(\R^n)$ directly yields $u_0 \in L^\infty(\R^n)$.
\end{remark}
\noindent\textbf{Proof of Proposition \ref{pro23}.} The proof can be divided into 2 steps. Step 1 investigates a semilinear parabolic equation and shows the local existence of the strong solution of Eqn. \er{star00}. Step 2 gives the uniqueness of the strong solution. \\

\noindent\textbf{Step 1 (Local existence).} In this step, we show the local existence of the strong solution, the proof is refined in spirit of \cite{BDP06,SK06}. Here, we denote $X_T$ by
\begin{align}
X_T:=\{f \in L^\infty(0,T;W^{2,n-1}(\R^n)), f_t \in L^{n-1}(Q_T)\big|~f \ge 0, \nonumber \\
\|f\|_{L^\infty(0,T;L^1 \cap L^\infty(\R^n))} \le C_1 \|u_0\|_{ L^\infty(0,T;L^1 \cap L^\infty(\R^n))}+C_2\}
\end{align}
for some $C_1,C_2$ are constants only depending on $n,\alpha,\beta,\sigma$ and $T>0$ to be determined later in Remark \ref{Tguji}. We also define
\begin{align}
W_u=\{u \in L^{n-1}(0,T;W^{2,n-1}(\R^n)) \mbox{~and~} u_t \in L^{n-1}(0,T;L^{n-1}(\R^n))\}.
\end{align}
Firstly, we consider
\begin{align}
V=K \ast f(x,t),\qquad x \in \R^n,0<t<T,
\end{align}
where $f \in X_T.$ Then by the weak Young inequality \cite{Lieb202}
\begin{align}
\nabla V \in L^\infty(0,T;L^{\frac{n}{n-1}}\cap L^\infty(\R^n)), \\
V \in L^\infty(0,T;L^{\frac{n}{n-2}}\cap L^\infty(\R^n)).
\end{align}
In addition, by the Maximum principle \cite{LSU} one has
\begin{align}
0 \le V(x,t) \le \|f\|_{L^\infty(Q_T)},~~~~x \in \R^n,~0<t<T.
\end{align}
Now we introduce the following problem
\begin{align}\label{E1}
\left\{
  \begin{array}{ll}
u_t=\Delta u-\sigma f^{\sigma-1} \nabla V\cdot \nabla u-f^{\sigma-1}\Delta V u+u^\alpha \left( 1-\int_{\R^n} f^\beta dx \right),\quad x \in \R^n,~0<t<T,\\
u\big|_{t=0}=u_0(x) \ge 0, \quad x \in \R^n.
  \end{array}
\right.
\end{align}
Let $u_0$ satisfies \er{initial}. Assume $f \in X_T,$ then $1-\int_{\R^n}f^\beta dx$ is bounded by a constant $C(\|u_0\|_{L^1 \cap L^\infty(\R^n)})$ which only depending on the initial data, then by virtue of \cite[Theorem 9.1]{LSU} and \cite{Fu66,H73,SK06}, equation \er{E1} corresponding to the initial data $u_0$ has a strong solution $u^{f} \in W_u$ and can be expressed by
\begin{align}
u(\cdot,t)&=G(\cdot,t)\ast u_0 - \sigma \int_0^t \left(f^{\sigma-1} \nabla V \cdot \nabla u(\cdot,s)\right)\ast G(\cdot,t-s) ds-\int_0^t \left(f^{\sigma-1} \Delta V u(\cdot,s)\right)\ast G(\cdot,t-s) ds \nonumber \\
&+\int_0^t \left( 1-\int_{\R^n} f^\beta dx \right) u^\alpha(\cdot,s)\ast G(\cdot,t-s) ds,
\end{align}
where $G(\cdot,t)$ is the Green function as in Lemma \ref{lem21}.

Next we define a mapping $\Phi$ by
\begin{align}
\Phi: f \in X_T \mapsto u^f \in W_u,
\end{align}
We show the solution $u^f$ is nonnegative as follows. Multiplying \er{E1} by $|u|^{k-2}u ~(k>1)$ and using Young's inequality we have
\begin{align}
&\frac{1}{k}\frac{d}{dt}\int_{\R^n} |u|^k dx+ (k-1) \int_{\R^n} |u|^{k-2} |\nabla u|^2dx \nonumber \\
=& -\int_{\R^n} \sigma f^{\sigma-1} \nabla V \cdot \nabla u |u|^{k-2}u dx-\int_{\R^n} f^{\sigma-1} \Delta V |u|^kdx+\int_{\R^n} u^{\alpha+1} |u|^{k-2} dx\left( 1-\int_{\R^n} f^\beta dx\right)\nonumber \\
\le &\frac{k-1}{4}\int_{\R^n} |u|^{k-2} |\nabla u|^2dx+\frac{\|\sigma f^{\sigma-1}\nabla V\|_{L^\infty(Q_T)}^2}{k-1} \int_{\R^n} |u|^k dx \nonumber \\
&+\|f^{\sigma-1}\Delta V\|_{L^\infty(Q_T)} \int_{\R^n} |u|^k dx +\left( 1+\|f\|_{L^\infty(0,T;L^1 \cap L^\infty(\R^n))}^\beta \right)\int_{\R^n} |u|^{\alpha+k-1} dx, \label{es1}
\end{align}
letting $F_0=1+\|f\|_{L^\infty(0,T;L^1 \cap L^\infty(\R^n))}^\beta$, we apply
\begin{align}
v=u^{k/2}, q=\frac{2(k+\alpha-1)}{k}, r=2, C_0=\frac{2(k-1)}{k^2}
\end{align}
in Lemma \ref{interpolation} for $k>\frac{n(\alpha-1)}{2}$ such that
\begin{align}
&F_0 \int_{\R^n} |u|^{\alpha+k-1} dx \nonumber\\
\le & \frac{2(k-1)}{k^2} \int_{\R^n} |\nabla u^{k/2}|^2dx+C\left(n,F_0 \right) \left( \frac{k^2}{2(k-1)} \right)^{\frac{\lambda q}{2-\lambda q}} \left( \int_{\R^n} |u|^k dx \right)^{\delta },
\end{align}
where
$$
\lambda=\frac{\frac{1}{2}-\frac{1}{q}}{\frac{1}{2}-\frac{1}{p}} \in (0,1),~\delta=\frac{1-\frac{q}{p}}{2-\frac{q}{2}-\frac{2}{p}}.
$$
Plugging it into \er{es1} one has
\begin{align}
&\frac{d}{dt} \|u\|_{L^k(\R^n)} \nonumber \\
\le &\left( \frac{\|\sigma f^{\sigma-1} \nabla V\|_{L^\infty(Q_T)}^2}{k-1} + \|f^{\sigma-1} \Delta V \|_{L^\infty(Q_T)} \right) \|u\|_{L^k(\R^n)}
+C\left(n,F_0 \right) \left( \frac{k^2}{2(k-1)} \right)^{C\left(\frac{1}{k}\right)} \|u\|_{L^k(\R^n)}^{k \delta-k+1},
\end{align}
taking $k \to \infty$ we derive
\begin{align}
\frac{d}{dt} \|u\|_{L^\infty(\R^n)} \le C\left(n,F_0,\|f^{\sigma-1}\Delta V\|_{L^\infty(Q_T)}\right) \|u\|_{L^\infty(\R^n)},
\end{align}
using Gronwall's inequality it's obtained that
\begin{align}\label{Linflocal}
\displaystyle \sup_{0<t<T} \|u(\cdot,t)\|_{L^\infty(\R^n)} \le \|u_0\|_{L^\infty(\R^n)}e^{C\left(n,F_0,\|f^{\sigma-1}\Delta V\|_{L^\infty(Q_T)}\right)T }.
\end{align}
Further integrating \er{E1} over $\R^n$ arrives at
\begin{align}\label{L1local}
\displaystyle \sup_{0<t<T} \|u(\cdot,t)\|_{L^1(\R^n)} \le \|u_0\|_{L^1(\R^n)} e^{C\left(\alpha,\|u\|_{L^\infty(Q_T)}\right)T}.
\end{align}
The nonnegativity of $u$ can be obtained by multiplying \er{E1} with $u^-:=-\min(u,0)$ that
\begin{align}
&\frac{1}{2}\frac{d}{dt} \int_{\R^n} |u^-|^2 dx \nonumber \\
&= -\int_{\R^n} |\nabla u^-|^2 dx-\sigma \int_{\R^n} f^{\sigma-1}\nabla V \cdot \nabla u^- u^- dx -\int_{\R^n} f^{\sigma-1} \Delta V |u^-|^2 dx +\left(1-\int_{\R^n} f^\beta dx \right)\int_{\R^n} u^\alpha u^- dx \nonumber \\
&\le -\frac{1}{2} \int_{\R^n} |\nabla u^-|^2 dx+ \frac{1}{2} \|\sigma f^{\sigma-1}\nabla V\|_{L^\infty(Q_T)}^2 \int_{\R^n}|u^-|^2 dx \nonumber \\
&+ \|f^{\sigma-1}\Delta V\|_{L^\infty(Q_T)} \int_{\R^n}|u^-|^2 dx+ F_0 \|u\|_{L^\infty(Q_T)}^{\alpha-1} \int_{\R^n} |u^-|^2 dx,
\end{align}
it follows
\begin{align}
\displaystyle \sup_{0<t<T} \|u(\cdot,t)\|_{L^2(\R^n)} \le e^{(\|\sigma f^{\sigma-1} \nabla V\|_{L^\infty(Q_T)}^2+2\|f^{\sigma-1}\Delta V\|_{L^\infty(Q_T)}+2F_0)T} \|u_0^-(\cdot,0)\|_{L^2(\R^n)}=0
\end{align}
directly assures that for all $0 \le t<T$
$$
u(x,t)\ge0,~~a.e.~~x \in \R^n.
$$
Furthermore, $\Phi$ is a contraction map in $L^\infty(0,T;L^{n-1}(\R^n))$. In fact, we consider the complete metric space $(X_T,d)$ where $d$ is defined by $d(f_1-f_2)=\|f_1-f_2\|_{L^\infty(0,T;L^{n-1}(\R^n))}$, we denote
\begin{align}
u_1=u^{f_1},~~u_2=u^{f_2},~~w=u_1-u_2,
\end{align}
from \er{E1} one has
\begin{align}\label{wstar}
&(u_1-u_2)_t=\Delta (u_1-u_2)-\left( \sigma f_1^{\sigma-1}\nabla V_1 \cdot \nabla u_1-\sigma f_2^{\sigma-1} \nabla V_2 \cdot \nabla u_2 \right) \nonumber \\
&-\left( f_1^{\sigma-1} \Delta V_1 u_1-f_2^{\sigma-1}\Delta V_2 u_2 \right)+\left(u_1^\alpha-u_2^\alpha\right) \left( 1-\int_{\R^n}f_2^\beta dx  \right)-u_1^\alpha \int_{\R^n} (f_1^\beta-f_2^\beta) dx,
\end{align}
the multiplication \er{wstar} by $|w|^{n-3}w$ and using H\"{o}lder's inequality give that for $\beta>n/2$
\begin{align}
&\frac{1}{n-1} \frac{d}{dt} \int_{\R^n} |w|^{n-1}dx +(n-2)\int_{\R^n} |\nabla w|^2 |w|^{n-3} dx \nonumber \\
&~\le \frac{n-2}{4} \int_{\R^n} |w|^{n-3}|\nabla w|^2 dx+\sigma \|f\|_{L^\infty(Q_T)}^{\sigma-1}\| u_1\|_{W^{2,n-1}(Q_T)}\int_{\R^n} |f_1-f_2|\cdot |w|^{n-2} dx \nonumber \\
&~~+\sigma \|f_2\|_{L^\infty(Q_T)}^\sigma \int_{\R^n} |\nabla w|\cdot |w|^{n-2}dx +\alpha \|u\|_{L^\infty(Q_T)}^{\alpha-1}F_0 \int_{\R^n} |w|^{n-1} dx+\int_{\R^n} \beta |f|^{\beta-1}|f_1-f_2| dx \int_{\R^n} u_1^\alpha |w|^{n-2} dx \nonumber \\
&~\le \frac{n-2}{4} \int_{\R^n} |w|^{n-3}|\nabla w|^2 dx+C\left(\|f\|_{L^\infty(Q_T)},\|u_1\|_{W^{2,n-1}(Q_T)}\right) \|f_1-f_2\|_{L^{n-1}(\R^n)} \| w^{n-2}\|_{L^{\frac{n-1}{n-2}}(\R^n)} \nonumber \\
&~~+\frac{n-2}{4} \int_{\R^n} |w|^{n-3}|\nabla w|^2 dx +C\left( \|f_2\|_{L^\infty(Q_T)},\|u\|_{L^\infty(Q_T)},F_0 \right) \int_{\R^n} |w|^{n-1} dx\nonumber \\
&~~+\beta \|f_1-f_2\|_{L^{n-1}(\R^n)} \|f^{\beta-1}\|_{L^{\frac{n-1}{n-2}}(\R^n)} \|w^{n-2}\|_{L^{\frac{n-1}{n-2}}(\R^n)} \|u_1^\alpha\|_{L^{n-1}(\R^n)}. \label{wcontraction}
\end{align}
Here $u$ and $f$ satisfy $\alpha u^{\alpha-1}=u_1^\alpha-u_2^\alpha$ and $\beta f^{\beta-1}=f_1^\beta-f_2^\beta$ by mean value theorem, \er{wcontraction} follows that
\begin{align}
&\frac{d}{dt} \|w\|_{L^{n-1}(\R^n)} \nonumber \\
\le & C_1\left( \|f\|_{L^\infty(Q_T)},\|u\|_{L^\infty(Q_T)}\right)\|w\|_{L^{n-1}(\R^n)}+C_2\left(
\|f\|_{L^\infty(Q_T)}, \|u_1\|_{W^{2,n-1}(Q_T)} \right)\|f_1-f_2\|_{L^{n-1}(\R^n)},
\end{align}
applying Gronwall's inequality yields
\begin{align}
\displaystyle \sup_{0<t<T} \|w(t)\|_{L^{n-1}(\R^n)} \le C_2 e^{C_1 T} \|f_1-f_2\|_{L^\infty(0,T;L^{n-1}(\R^n))},
\end{align}
hence there exists $T_\ast=T_\ast\left(\|u_0\|_{L^1 \cap W^{2,n-1}(\R^n)},\|u_0\|_{L^\infty(Q_T)},T\right)\le T$ small such that
\begin{align}
\displaystyle \sup_{0<t<T_\ast} \|w(t)\|_{L^{n-1}(\R^n)} \le \frac{1}{2} \|f_1-f_2\|_{L^\infty(0,T_\ast;L^{n-1}(\R^n))}.
\end{align}
Using Banach's fixed point theorem, we have that $\Phi$ has a fixed point $\Phi(f)=u^f=f \in X_{T_\ast}$. Iterating the method we prove the existence of the strong solution $u$ of equation on $[0,T)$
\begin{align}\label{232}
\left\{
  \begin{array}{ll}
 u_t=\Delta u-\sigma u^{\sigma-1} \nabla V \cdot \nabla u-u^{\sigma-1} \Delta V u+u^\alpha \left( 1-\int_{\R^n} u^\beta dx \right),\quad   &x \in \R^n, 0<t<T,
\\
v=K \ast u, &x \in \R^n, 0<t<T,
\\
 u(x,0)=u_0(x)\geq 0,\quad & x \in \R^n.
  \end{array}
\right.
\end{align}
\\

\noindent\textbf{Step 2 (Uniqueness).} The uniqueness of the strong solution can be shown in the follows. Assume $(u_1,v_1)$ and $(u_2,v_2)$ solve \er{232} in $Q_{T}$ with initial data \er{initial}, then
\begin{align}
(u_1-u_2)_t=\Delta(u_1-u_2)-\nabla \cdot \left((u_1^\sigma-u_2^\sigma) \nabla v_1+u_2^\sigma \nabla(v_1-v_2)\right) \nonumber\\
+u_1^\alpha-u_2^\alpha+\left(u_2^\alpha-u_1^\alpha\right)\int_{\R^n} u_1^\beta dx+u_2^\alpha \int_{\R^n} \left(u_2^\beta-u_1^\beta\right) dx, \label{233}
\end{align}
multiplying \er{233} with $u_1-u_2$ we obtain
\begin{align}\label{I1}
&\frac{1}{2}\frac{d}{dt}\int_{\R^n} |u_1-u_2|^2 dx \nonumber \\
&=-\int_{\R^n} |\nabla (u_1-u_2)|^2dx+\int_{\R^n} (u_1^\sigma-u_2^\sigma) \nabla(u_1-u_2)\cdot \nabla v_1 dx+\int_{\R^n} u_2^\sigma \nabla(u_1-u_2)\cdot\nabla(v_1-v_2) dx\nonumber \\
&~~+\left(1-\int_{\R^n}u_1^\beta dx\right)\int_{\R^n} (u_1^\alpha-u_2^\alpha)(u_1-u_2)dx+ \int_{\R^n} \left(u_2^\beta-u_1^\beta\right) dx \int_{\R^n} u_2^\alpha (u_1-u_2) dx \nonumber \\
&\le -\int_{\R^n} |\nabla(u_1-u_2)|^2 dx+\frac{1}{4}\int_{\R^n} |\nabla(u_1-u_2)|^2 dx+\|\nabla v_1\|_{L^\infty(Q_{T})}^2 \left( \sigma \|u\|_{L^\infty(Q_{T})}^{\sigma-1}\right)^2 \int_{\R^n} |u_1-u_2|^2dx \nonumber \\
&~~+\frac{1}{4}\int_{\R^n} |\nabla(u_1-u_2)|^2 dx+\int_{\R^n}|u_2^\sigma \nabla(v_1-v_2)|^2 dx+\alpha \left(1+\|u_1\|_{L^\infty(0,T;L^\beta(\R^n))}^\beta \right)\|u\|_{L^\infty(Q_{T})}^{\alpha-1} \int_{\R^n} |(u_1-u_2)|^2 dx\nonumber \\ &~~+\int_{\R^n}|u_2^\beta-u_1^\beta|dx \int_{\R^n}|u_2^\alpha(u_1-u_2)|dx.
\end{align}
Here $u$ comes from mean value theorem by $\alpha u^{\alpha-1}=u_1^\alpha-u_2^\alpha$. Therefore
\begin{align}
&\frac{d}{dt}\int_{\R^n} |u_1-u_2|^2 dx \le C\big(\|\nabla v_1\|_{L^\infty(Q_{T})},\|u\|_{L^\infty(Q_{T})}, \|u_1\|_{L^\infty(0,T;L^\beta(\R^n))} \big)\int_{\R^n} |u_1-u_2|^2dx \nonumber \\
&+ \int_{\R^n}|u_2^\sigma \nabla(v_1-v_2)|^2dx+\int_{\R^n}|u_2^\beta-u_1^\beta|dx \int_{\R^n}|u_2^\alpha(u_1-u_2)|dx
=I_1+I_2+I_3,
\end{align}
By H\"{o}lder's inequality and weak Young's inequality \cite{Lieb202} we have
\begin{align}\label{I2}
I_2 &\le \big\||\nabla(v_1-v_2)|^2\big\|_{L^{\frac{n}{n-2}}(\R^n)}\|u_2^{2\sigma}\|_{L^{\frac{n}{2}}(\R^n)} \nonumber\\
&\le C \big\|\frac{x-y}{|x-y|^n} \big\|^2_{L_w^{\frac{n}{n-1}}(\R^n)} \|u_1-u_2\|_{L^2(\R^n)}^2 \|u_2^{2 \sigma}\|_{L^{\frac{n}{2}}(\R^n)},
\end{align}
and
\begin{align}\label{I3}
I_3 &\le \int_{\R^n} \beta u^{\beta-1}|u_1-u_2|dx \int_{\R^n} u_2^\alpha |u_1-u_2|dx \nonumber \\
&\le \|u_1-u_2\|_{L^2(\R^n)} \|\beta u^{\beta-1}\|_{L^2(\R^n)} \|u_2^\alpha \|_{L^2(\R^n)} \|u_1-u_2\|_{L^2(\R^n)} \nonumber \\
&\le C\left(\|u\|_{L^1\cap L^\infty(\R^n)}\right) \|u_1-u_2\|_{L^2(\R^n)}^2,  ~~(\beta \ge 3/2)
\end{align}
taking \er{I1}, \er{I2} and \er{I3} together one has
\begin{align}
\frac{d}{dt}\int_{\R^n} |u_1-u_2|^2 dx\le C(\|u\|_{L^1\cap L^\infty(\R^n)},n) \int_{\R^n} |u_1-u_2|^2 dx,
\end{align}
this yields $u_1=u_2$ in $Q_{T}$ which implies the uniqueness of solutions. Thus we complete the proof of local existence and uniqueness of the strong solution. $\Box$

\begin{remark}\label{Tguji}
In Proposition \ref{pro23}, the bounded time $T$ in $X_T$ can be preestimated as follows
\begin{align}
f_t&=\Delta f-\nabla \cdot(f^\sigma \nabla v)+ f^\alpha \left(1-\int_{\R^n}f^\beta dx\right), \label{333} \\
v&=K \ast f, \nonumber \\
f(0)&=u_0 \ge 0. \nonumber
\end{align}
multiplying \er{333} by $rf^{r-1}(r>1)$ obtains that
\begin{align}
\frac{d}{dt}\int_{\R^n} f^r dx &=-\frac{4(r-1)}{r} \int_{\R^n}|\nabla f^{\frac{r}{2}}|^2 dx +\frac{r(r-1)}{\sigma+r-1}\int_{\R^n} f^{\sigma+r} dx+r (1-\int_{\R^n} f^\beta dx) \int_{\R^n} f^{\alpha+r-1} dx \nonumber \\
&\le \frac{r(r-1)}{\sigma+r-1} \|f\|_{L^\infty(\R^n)}^\sigma \|f\|_{L^r(\R^n)}^r+r \|f\|_{L^\infty(\R^n)}^{\alpha-1}\|f\|_{L^r(\R^n)}^r,
\end{align}
thus
\begin{align}
\frac{d}{dt} \|f\|_{L^r(\R^n)} \le \frac{r-1}{\sigma+r-1}\|f\|_{L^\infty(\R^n)}^\sigma \|f\|_{L^r(\R^n)}+ \|f\|_{L^\infty(\R^n)}^{\alpha-1}\|f\|_{L^r(\R^n)},
\end{align}
letting $r \to \infty$ one has
\begin{align}
\|f\|_{L^\infty(\R^n)} \le \|u_0\|_{L^\infty(\R^n)}+\int_0^t \|f\|_{L^\infty(\R^n)}^{\sigma+1} ds+\int_0^t \|f\|_{L^\infty(\R^n)}^\alpha ds,
\end{align}
hence from ODE inequality we have that there is a maximum existence time $T=T(\|u_0\|_{L^\infty(\R^n)})$ such that $f$ is bounded from above in $[0,T)$.
\end{remark}

\section{Proof of Theorem \ref{thm1}}\label{sec4}
\def\theequation{4.\arabic{equation}}\makeatother
\setcounter{equation}{0}
\def\thetheorem{4.\arabic{thm}}\makeatother

In this section, we derive the a priori estimates of the strong solution and complete the proof of Theorem \ref{thm1}.
\begin{proposition}\label{pro25}
Let $\sigma \ge 1, \beta>1, \alpha>1,$ $u_0$ satisfies \er{initial}. If
\begin{align}
\sigma+1\le \alpha \mbox{ and } \alpha<1+2\beta/n
\end{align}
or
\begin{align}
\sigma+1>\alpha \mbox{ and } (\sigma+1-\alpha)(n+2)<2\beta-n(\alpha-1),
\end{align}
then for any $0<t<T_{\max}$, the solution of \er{star00} satisfies that \\
\noindent(1) For $\beta+\alpha-1 \le k<\infty,$
\begin{align}
\|u(\cdot,t)\|_{L^k(\R^n)} \le C(\|u_0\|_{L^k(\R^n)},k).
\end{align}
\noindent(2) The uniformly boundedness of solution $$\|u(\cdot,t)\|_{L^\infty(\R^n)}\le C \left(\|u_0\|_{L^{\beta+\alpha-1}(\R^n)},\|u_0\|_{L^\infty(\R^n)}\right),$$ where $C$ is a positive constant depending on $\|u_0\|_{L^{\beta+\alpha-1}(\R^n)}$ and $\|u_0\|_{L^\infty(\R^n)}$ but not on $T_{\max}$. \\
\noindent(3) For $1 \le k<\beta+\alpha-1$,
\begin{align}
\|u(\cdot,t)\|_{L^k(\R^n)} \le C(\|u_0\|_{L^k(\R^n)},T_{\max})
\end{align}
where $C$ is a positive constant depending on $\|u_0\|_{L^k(\R^n)}$ and $T_{\max}$. \\

Especially, when $\sigma+1=\alpha,$ we have that for $\beta \le k \le \infty$
\begin{align}
\|u(\cdot,t)\|_{L^k(\R^n)} \le C(\|u_0\|_{L^k(\R^n)}),
\end{align}
where $C$ only depends on $\|u_0\|_{L^k(\R^n)}$ not on $T_{\max}$.
\end{proposition}
\noindent{\bf Proof of Proposition \ref{pro25}.} For the rigorous proof, we should multiply \er{star00} by $ku^{k-1}\psi_l$, where $\psi_l$ is a standard cut-off function. By the limiting process we can justify the following formal calculation. Throughout the proof, we suppose $\sigma+1=\eta.$ \\

\noindent\textbf{Step 1 (A priori estimates).} Firstly multiplying \er{star00} with $k u^{k-1}(k \ge 1)$ one has
\begin{align}\label{baseguji}
& \frac{d}{dt} \int_{\R^n}  u^k dx  + \frac{4(k-1)}{k}
\int_{\R^n} |\nabla u^{\frac{k}{2}} |^2 dx + k \int_{\R^n} u^{k+\alpha-1} dx \int_{\R^n} u^\beta dx \nonumber\\
=& ~k \int_{\R^n} u^{k+\alpha-1} dx +\frac{k(k-1)}{k+\sigma-1} \int_{\R^n} u^{\eta+k-1} dx.
\end{align}
In order to control the right hand side of \er{baseguji} by using the two nonnegative terms in the left hand side of \er{baseguji}, we apply
$$~v=u^{k/2}, q=\frac{2(k+\alpha-1)}{k}, r=\frac{2k'}{k}, C_0=\frac{k-1}{k^2} $$
in Lemma \ref{interpolation} with
$k>\frac{2(\alpha-1)}{p-2}$ (which is $q<p$) and $\frac{p(\alpha-1)}{p-2}<k'<k+\alpha-1$ (which is $\frac{q}{r}<\frac{2}{r}+1-\frac{2}{p},r<q$)
such that
\begin{align}\label{1}
\int_{\R^n} u^{\alpha+k-1} dx \le \frac{k-1}{k^2} \| \nabla u^{\frac{k}{2}} \|_{L^2(\R^n)}^2+C(k) \|u\|_{L^{k'}(\R^n)}^{b_\alpha},
\end{align}
where
$$ b_\alpha=\frac{(1-\lambda_\alpha)(\alpha+k-1)}{1-\frac{\lambda_\alpha(\alpha+k-1)}{k}},~\lambda_\alpha=\frac{\frac{k}{2k'}-\frac{k}{2(\alpha+k-1)}}{\frac{k}{2k'}-\frac{1}{p}}\in(0,1).
$$
Similarly, taking
$$~v=u^{k/2}, q=\frac{2(k+\eta-1)}{k}, r=\frac{2k'}{k}, C_0=\frac{k+\sigma-1}{k^2} $$
in Lemma \ref{interpolation} with $k>\frac{2(\eta-1)}{p-2}$ (which is $q<p$) and $\frac{p(\eta-1)}{p-2}<k'<k+\eta-1$ (which is $\frac{q}{r}<\frac{2}{r}+1-\frac{2}{p},r<q$) leads to
\begin{align}\label{2}
\int_{\R^n} u^{\eta+k-1} dx \le \frac{k+\sigma-1}{k^2} \| \nabla u^{\frac{k}{2}} \|_{L^2(\R^n)}^2+C(k) \|u\|_{L^{k'}(\R^n)}^{b_\eta},
\end{align}
where
$$ b_\eta=\frac{(1-\lambda_\eta)(\eta+k-1)}{1-\frac{\lambda_\eta(\eta+k-1)}{k}},~\lambda_\eta=\frac{\frac{k}{2k'}-\frac{k}{2(\eta+k-1)}}{\frac{k}{2k'}-\frac{1}{p}} \in(0,1).$$
Hence taking \er{1} and \er{2} together we will conduct further estimates of \er{baseguji} for
\begin{align}\label{k}
k>\max \left( \frac{2(\eta-1)}{p-2},\frac{2(\alpha-1)}{p-2},1 \right)
\end{align}
and
\begin{align}\label{kprime}
\max \left( \frac{p(\eta-1)}{p-2},\frac{p(\alpha-1)}{p-2},1,\frac{k}{2} \right)<k'<\min\left( k+\alpha-1, k+\eta-1 \right)
\end{align}
where we have used the fact that $1 \le r$ is equivalent to $\frac{k}{2} \le k'.$

Combining \er{1} and \er{2} we infer from \er{baseguji} that
\begin{align}\label{guji1}
& \frac{d}{dt}\int_{R^n} u^k dx+ k \int_{\R^n} u^\beta dx \int_{\R^n} u^{k+\alpha-1} dx+\frac{2(k-1)}{k} \|\nabla u^{\frac{k}{2}} \|_{L^2(\R^n)}^2 \nonumber \\
\le ~& C(k,\alpha)  \|u\|_{L^{k'}(\R^n)}^{b_\alpha}+ C(k,\eta) \|u\|_{L^{k'}(\R^n)}^{b_\eta}.
\end{align}
We further assume $\beta<k'$ and use the following interpolation inequalities such that
\begin{align}\label{guji2}
\|u\|_{L^{k'}}^{b_\alpha} \le \left( \|u\|_{L^{k+\alpha-1}}^{k+\alpha-1} \|u\|_{L^\beta}^\beta \right)^{\frac{b_\alpha \theta}{k+\alpha-1}} \|u\|_{L^\beta}^{b_\alpha(1-\theta-\frac{\theta \beta}{k+\alpha-1})}
\end{align}
and
\begin{align}\label{guji3}
\|u\|_{L^{k'}}^{b_\eta} \le \left( \|u\|_{L^{k+\alpha-1}}^{k+\alpha-1} \|u\|_{L^\beta}^\beta \right)^{\frac{b_\eta \theta}{k+\alpha-1}} \|u\|_{L^\beta}^{b_\eta(1-\theta-\frac{\theta \beta}{k+\alpha-1})}
\end{align}
where
$$
\theta=\frac{\frac{1}{\beta}-\frac{1}{k'}}{\frac{1}{\beta}-\frac{1}{k+\alpha-1}}.
$$
To use Young's inequality, we need the following three conditions that
\begin{align}\label{star10}
1-\theta-\frac{\theta \beta}{k+\alpha-1}=0
\end{align}
and
\begin{align} \label{star}
\frac{b_\alpha \theta }{k+\alpha-1}<1
\end{align}
as well as
\begin{align}\label{starstar}
\frac{b_\eta \theta }{k+\alpha-1}<1.
\end{align}
Firstly for \er{star10}, thanks to the arbitrariness of $k',$ we take
$$
k'=\frac{k+\alpha-1+\beta}{2} \in (\beta,k+\alpha-1)
$$
so that \er{star10} holds true. Here combining \er{k} and \er{kprime}, $k$ satisfies
\begin{align}\label{kxin}
k>\max \left( \frac{2(\eta-1)}{p-2},~\frac{2(\alpha-1)}{p-2}, ~\frac{2p(\eta-1)}{p-2}-\beta-(\alpha-1),~ \frac{2p(\alpha-1)}{p-2}-\beta-(\alpha-1)\right)=:K_0.
\end{align}
Next recalling $b_{\alpha}$ and $\theta$, \er{star} is equivalent to
\begin{align}
(1-\lambda_{\alpha})\left( \frac{1}{\beta}-\frac{1}{k'}  \right)< \left( 1-\frac{\lambda_\alpha (k+\alpha-1)}{k} \right)\left( \frac{1}{\beta}-\frac{1}{k+\alpha-1}\right),
\end{align}
it also reads
\begin{align}
\lambda_\alpha-1<\frac{\lambda_\alpha}{k} k'-\frac{k'}{k+\alpha-1}-\frac{\lambda_\alpha(\alpha-1)}{k \beta} k',
\end{align}
substituting $\lambda_\alpha$ into the above formula one has that for
\begin{align}\label{tiaojian1}
1 \le \alpha<1+\left( 1-\frac{2}{p} \right) \beta,
\end{align}
\er{star} holds true for any $k$ satisfying \er{kxin}.
For \er{starstar}, under the condition \er{tiaojian1} similar arguments arrive at
\begin{align}
\left( \frac{k}{2}-\frac{k+\eta-1}{p} \right) \left( \frac{k'}{\beta}-1  \right)<\left( \left(\frac{1}{2}-\frac{1}{p} \right)k'-\frac{\eta-1}{2}  \right) \left( \frac{\alpha-1}{\beta}-1+\frac{k}{\beta} \right),
\end{align}
it can be written as
\begin{align}\label{starstarstar}
\frac{\eta-\alpha}{\beta} ~\left( \frac{k}{2}-\frac{k'}{p} \right)<(\eta+k-1-k') ~\left( \frac{1}{2}-\frac{1}{p}-\frac{\alpha-1}{2 \beta} \right).
\end{align}

Denote
$$
A_0=\frac{1}{2}-\frac{1}{p}-\frac{\alpha-1}{2 \beta}>0,~A_1=\frac{\eta-\alpha}{\beta},
$$
recalling \er{k} and \er{kprime}, if $\eta \le \alpha,$ then $A_0>0$ is enough to guarantee that \er{starstar} holds true. Otherwise if $\eta > \alpha,$ plugging
$
k'=\frac{k+\alpha-1+\beta}{2}
$
into \er{starstarstar} yields
$$
A_0(\eta-1)+\frac{A_1}{2p}(\alpha+\beta-1)+\Big[  \frac{A_0}{2}-A_1 \left( \frac{1}{2}-\frac{1}{2p} \right)  \Big] k~>~A_0 \frac{\alpha+\beta-1}{2}.
$$
Hence if
$$
\frac{A_0}{2}-A_1 \left( \frac{1}{2}-\frac{1}{2p} \right)>0,
$$
or equivalently
\begin{align}\label{tiaojian2}
\frac{\eta-\alpha}{\beta} \left(1-\frac{1}{p} \right)< \frac{1}{2}-\frac{1}{p}-\frac{\alpha-1}{2 \beta},
\end{align}
then \er{starstar} holds true for
\begin{align}\label{DB}
k>\frac{A_0 \left( \frac{\alpha+\beta-1}{2}-(\eta-1) \right) -\frac{A_1}{2p}(\alpha+\beta-1)}{\frac{A_0}{2}-A_1\left( \frac{1}{2}-\frac{1}{2p} \right)}=:\frac{D}{B}.
\end{align}
In the following, we prove $D>0$. In fact,
\begin{align}
D&=A_0 \left( \frac{\alpha+\beta-1}{2}-(\eta-1) \right) -\frac{A_1}{2p}(\alpha+\beta-1) \nonumber\\
& > A_0 \left(\frac{\alpha+\beta-1}{2}-\frac{p-2}{p-1} \frac{\alpha+\beta-1}{2} \right)-\frac{A_1}{2p}(\alpha+\beta-1) \nonumber\\
&= \frac{A_0}{p-1}\frac{\alpha+\beta-1}{2} -\frac{A_1}{2p}(\alpha+\beta-1) \nonumber \\
&>0,
\end{align}
where we have used \er{tiaojian2} and its reformulation
\begin{align}\label{tiaojian2bianxing}
\eta-1<\frac{p-2}{2(p-1)} (\alpha-1+\beta).
\end{align}
Furthermore, after some calculations we have
\begin{align}
\frac{D}{B} & = \frac{\left(-A_0-\frac{\alpha+\beta-1}{2p\beta}  \right) \eta + \frac{\alpha}{p \beta}\frac{\alpha+\beta-1}{2}+A_0\frac{\alpha+\beta+1}{2} }{ -\frac{1}{\beta} \left( \frac{1}{2}-\frac{1}{2p} \right) \eta+\frac{A_0}{2}+\frac{\alpha}{\beta} \left( \frac{1}{2}-\frac{1}{2p} \right)   }   \nonumber \\[2mm]
& = \frac{-\frac{1}{\beta}(\beta-(\alpha-1)) \left( \frac{1}{2}-\frac{1}{2p} \right) \eta+ \big(\beta-(\alpha-1)\big) \left(\frac{A_0}{2}+\frac{\alpha}{\beta} \left( \frac{1}{2}-\frac{1}{2p} \right) \right)}{ -\frac{1}{\beta} \left( \frac{1}{2}-\frac{1}{2p} \right) \eta+\frac{A_0}{2}+\frac{\alpha}{\beta} \left( \frac{1}{2}-\frac{1}{2p} \right)   } \nonumber\\
&=\beta-(\alpha-1).
\end{align}

Now coming back to \er{guji2} and \er{guji3}, we infer from \er{star} and \er{starstar} by using Young's inequality that
\begin{align}
C(k,\alpha)\|u\|_{L^{k'}}^{b_\alpha} & \le \left( \|u\|_{L^{k+\alpha-1}}^{k+\alpha-1} \|u\|_{L^\beta}^\beta \right)^{\frac{b_\alpha \theta}{k+\alpha-1}} \nonumber\\
& \le \frac{k}{4} \|u\|_{L^{k+\alpha-1}(\R^n)}^{k+\alpha-1} \|u\|_{L^\beta(\R^n)}^\beta+C_1(k)
\end{align}
and
\begin{align}
C(k,\eta) \|u\|_{L^{k'}}^{b_\eta} & \le \left( \|u\|_{L^{k+\alpha-1}}^{k+\alpha-1} \|u\|_{L^\beta}^\beta \right)^{\frac{b_\eta \theta}{k+\alpha-1}} \nonumber \\
& \le \frac{k}{4} \|u\|_{L^{k+\alpha-1}(\R^n)}^{k+\alpha-1} \|u\|_{L^\beta(\R^n)}^\beta+C_2(k).
\end{align}

Therefore, substituting the above arguments into \er{guji1} we thus end up with for all $0<t<T_{\max}$
\begin{align}\label{guji}
\frac{d}{dt}\int_{\R^n} u^k dx+\frac{k}{2} \int_{\R^n} u^\beta dx \int_{\R^n} u^{k+\alpha-1} dx+\frac{2(k-1)}{k} \|\nabla u^{\frac{k}{2}} \|_{L^2(\R^n)}^2 \le C(k)
\end{align}
for $k>\beta-(\alpha-1)$ from \er{kxin}, \er{DB} and \er{tiaojian1}, \er{tiaojian2}. Precisely, for $\eta>\alpha$, \er{tiaojian2} is equivalent to
\begin{align}
\frac{2(\eta-1)}{p-2}(p-1)<\beta+\alpha-1,
\end{align}
and then $k>\max\{K_0,\beta-(\alpha-1)\}=\beta-(\alpha-1).$ Otherwise for $\eta \le \alpha,$ by virtue of \er{tiaojian1} it leads to
$k>\max\{ K_0, \beta-(\alpha-1) \}=\beta-(\alpha-1)$. \\

\noindent\textbf{Step 2 ($L^k(\R^n)$ estimates for $\beta+\alpha-1 \le k \le \infty$).} Firstly we note that for $\eta \le \alpha$, then $\beta+\alpha-1>\frac{2(\alpha-1)}{p-2}$ with the help of \er{tiaojian1}, therefore we can take $k=\beta+\alpha-1$ in \er{guji} and by H\"{o}lder inequality and Young's inequality on has
\begin{align}
\|u\|_{L^{\beta+\alpha-1}(\R^n)}^{\beta+\alpha-1} & \le \left( \|u\|_{L^{\beta+2(\alpha-1)}(\R^n)}^{\beta+2(\alpha-1)} \|u\|_{L^\beta(\R^n)}^\beta \right)^{1/2} \nonumber \\
& \le \frac{\beta+\alpha-1}{2} \|u\|_{L^{\beta+2(\alpha-1)}(\R^n)}^{\beta+2(\alpha-1)} \|u\|_{L^\beta(\R^n)}^\beta +\frac{1}{2(\beta+\alpha-1)},
\end{align}
plugging the above formula into \er{guji} we have
\begin{align}
\frac{d}{dt} \int_{\R^n} u^{\beta+\alpha-1} dx+ \int_{\R^n} u^{\beta+\alpha-1} \le C(\beta+\alpha-1)
\end{align}
which follows the uniformly boundedness in time
\begin{align}
\int_{\R^n} u^{\beta+\alpha-1} dx \le \|u_0\|_{L^{\beta+\alpha-1}(\R^n)}^{\beta+\alpha-1} +C(\alpha,\beta).
\end{align}
On the other hand, letting
\begin{align}
v=u^{k/2},~q=2,~1\le r<2,~C_0=\frac{2(k-1)}{k}
\end{align}
in Lemma \ref{interpolation} one has that for $n \ge 1$
\begin{align}\label{336}
\int_{\R^n} u^k dx \le \frac{2(k-1)}{k} \| \nabla u^{\frac{k}{2}} \|_{L^2(\R^n)}^2 +C(n,k) \|u\|_{L^{k_1}(\R^n)}^k
\end{align}
where $k_1=\frac{kr}{2} <k$. Furthermore, for $\beta<k_1<k+\alpha-1$ we can take $k_1=\frac{\beta+\alpha-1+k}{2}<k$ which is
\begin{align}
k>\beta+\alpha-1
\end{align}
so that
\begin{align}\label{339}
\|u\|_{L^{\frac{\beta+\alpha-1+k}{2}}(\R^n)}^k \le \left( \|u\|_{L^{k+\alpha-1}(\R^n)}^{k+\alpha-1} \|u\|_{L^\beta(\R^n)}^\beta \right)^{\frac{k}{\beta+\alpha-1+k}}.
\end{align}
Combining \er{336} and \er{339} together yields
\begin{align}\label{340}
\int_{\R^n} u^k dx &\le \frac{2(k-1)}{k} \| \nabla u^{\frac{k}{2}} \|_{L^2(\R^n)}^2+C(n,k) \left( \|u\|_{L^{k+\alpha-1}(\R^n)}^{k+\alpha-1} \|u\|_{L^\beta(\R^n)}^\beta \right)^{\frac{k}{\beta+\alpha-1+k}} \nonumber \\
& \le \frac{2(k-1)}{k} \| \nabla u^{\frac{k}{2}} \|_{L^2(\R^n)}^2+\frac{k}{2}  \|u\|_{L^{k+\alpha-1}(\R^n)}^{k+\alpha-1} \|u\|_{L^\beta(\R^n)}^\beta +C(n,k).
\end{align}
Substituting \er{340} into \er{guji} one has
\begin{align}
\frac{d}{dt} \int_{\R^n} u^k dx + \int_{\R^n} u^k dx \le C(n,k).
\end{align}
It can obtained that for any $\beta+\alpha-1<k<\infty$
\begin{align}
\int_{\R^n} u^k dx \le C(\|u_0\|_{L^k(\R^n)},k).
\end{align}

Furthermore, for the $L^\infty$ norm, according to \er{baseguji},
we can conduct similar procedures as Step 4 in \cite{BL16} in terms of $\int_{\R^n} u^{k+\alpha-1} dx$ and $\int_{\R^n} u^{k+\eta-1} dx$ together on the right hand side of \er{baseguji} and get
\begin{align}\label{inf}
\|u\|_{L^\infty(\R^n)} \le C(\alpha,\eta,\|u_0\|_{L^\beta(\R^n)},\|u_0\|_{L^\infty(\R^n)}).
\end{align}
\\

\noindent\textbf{Step 3 ($L^\beta$ estimates).}
When $\sigma+1=\alpha,$ letting $k=\beta$ in \er{baseguji} one has
\begin{align}
\frac{d}{dt}\int_{\R^n} u^\beta dx \le \beta \int_{\R^n} u^{\beta+\alpha-1} dx \left(1+\frac{\beta-1}{\beta+\sigma-1}-\int_{\R^n} u^\beta dx \right),
\end{align}
it follows that for all $0<t<T_{\max}$
\begin{align}
\int_{\R^n} u(\cdot,t)^\beta dx \le \max \left(\int_{\R^n}u_0^\beta dx,1+\frac{\beta-1}{\beta+\sigma-1}\right).
\end{align}
\\

\noindent\textbf{Step 4 ($L^k$ estimates for $1 \le k<\beta+\alpha-1$).} By virtue of \er{guji}, we have that for any $0<t<T_{\max}$ and $\beta-(\alpha-1)<k<\beta+\alpha-1$,
\begin{align}\label{beal}
\int_{\R^n} u(\cdot,t)^k dx\le C(k)T_{\max}+\int_{\R^n} u_0^k dx.
\end{align}
Furthermore, integrating \er{star00} over $\R^n$ and using \er{inf} and \er{beal} get
\begin{align}
\frac{d}{dt} \int_{\R^n} u dx &= \int_{\R^n} u^\alpha dx \left( 1-\int_{\R^n} u^\beta dx  \right) \nonumber \\
& \le \|u\|_{L^\infty(\R^n)}^{\alpha-1}\left(1+\|u\|_{L^\beta(\R^n)}\right) \int_{\R^n} u dx \nonumber \\
& \le C\left(\|u_0\|_{L^\infty(\R^n)},\|u_0\|_{L^\beta(\R^n)},T_{\max} \right) \int_{\R^n} u dx,
\end{align}
thus for any $0<t<T_{\max}$
\begin{align}
\|u(\cdot,t)\|_{L^1(\R^n)}\le e^{C\left(\|u_0\|_{L^\infty(\R^n)},\|u_0\|_{L^\beta(\R^n)},T_{\max} \right) T_{\max}} \|u_0\|_{L^1(\R^n)}.
\end{align}
Consequently we conclude that for $1\le k<\beta+\alpha-1$ and any $0<t<T_{\max}$
\begin{align}
\|u(\cdot,t)\|_{L^k(\R^n)} \le C(\|u_0\|_{L^k(\R^n)},\beta,\alpha,T_{\max}).
\end{align}
This completes the a priori estimates. $\Box$ \\

Proposition \ref{pro25} together with the blow-up criterion \er{blowup} allow us to state without further arguments
\noindent{\bf Proof of Theorem \ref{thm1}.} With the aid of the blow-up criterion \er{blowup} and the uniformly boundedness of the solution in Proposition \ref{pro25}, there exists a positive constant $C\left(\|u_0\|_{L^{\beta+\alpha-1}(\R^n)},\|u_0\|_{L^\infty(\R^n)}\right)$ such that
\begin{align}
\|u(\cdot,t)\|_{L^\infty(\R^n)} \le C~~\mbox{ for all } 0<t<\infty.
\end{align}
By Proposition \ref{pro23} we obtain the desired result. The proof of Theorem \ref{thm1} is completed. $\Box$

\section{Conclusions}\label{sec5}
\def\theequation{5.\arabic{equation}}\makeatother
\setcounter{equation}{0}
\def\thetheorem{5.\arabic{thm}}\makeatother
\setcounter{thm}{0}

This paper concerns Eq. \er{star00} in terms of different reaction and aggregation exponents. If $\alpha \ge \sigma+1,$ then the growth in reaction dominates and letting $\alpha<1+2\beta/n$ in the death $u^\alpha \int_{\R^n} u^\beta dx$ can prevent blow-up. While for $\alpha < \sigma+1,$ the aggregation dominates and let the aggregation exponent $\sigma+1<(2\alpha+2\beta+n)/(n+2)$ thus the solution will exist globally. Moreover, if $\sigma+1=\alpha,~n(\alpha-1)/\beta=2,$ then $u_\lambda(x,t)=\lambda^{\frac{n}{\beta}}u(\lambda x, \lambda^2 t)$ is also a solution of \er{star00} and the scaling preserves the $L^\beta$ norm in space. When $\frac{n}{\beta}(\alpha-1)<2,$ for low density (small $\lambda$), the aggregation dominates the diffusion thus prevents spreading. While for high density (large $\lambda$), the diffusion dominates the aggregation and thus blow-up is precluded. Hence, in this case, the solution will exist globally (Theorem \ref{thm1}) and we believe that the solution converges to the stationary solution as time goes to infinity. On the contrary, both global existence and finite time blow-up may occur for $\frac{n}{\beta}(\alpha-1)>2,$ hence our conjecture is that there exists finite time blow-up for $\alpha-1>2\beta/n.$ As to the case $\alpha-1=2\beta/n$, similar to \cite{Jose2009}, whether there is a critical value for the initial data sharply separating global existence and finite time blow-up is also unknown. Our result is to be considered as the first step to a more general theory of chemotaxis system with nonlocal nonlinear reaction. This will be a fertile area to explore.

\end{document}